\documentclass[preprint,12pt]{elsarticle}

\usepackage{amsmath,amssymb,amsthm}
\usepackage{float}
\usepackage{hyperref}
\usepackage{xcolor}
\usepackage{subcaption}
\usepackage{algorithm}
\usepackage{algpseudocode}
\usepackage{algorithmicx}
\usepackage{placeins}

\algdef{SE}[DOWHILE]{Do}{doWhile}{\algorithmicdo}[1]{\algorithmicwhile\ #1}
\algrenewcommand\algorithmicrequire{\textbf{Input:}}
\algrenewcommand\algorithmicensure{\textbf{Output:}}
\makeatletter
\def\BState{\State\hskip-\ALG@thistlm}
\makeatother

\theoremstyle{plain}
\newtheorem{thm}{Theorem}[section]
\newtheorem{lem}[thm]{Lemma}

\newtheorem{prop}[thm]{Proposition}

\theoremstyle{definition}
\newtheorem{defn}[thm]{Definition}

\theoremstyle{remark}
\newtheorem{rem}[thm]{Remark}

\numberwithin{equation}{section}

\def\ff{{\mathcal F}}
\newcommand{\g}[1]{\textbf{#1}}
\newcommand{\norm}[1]{{\left\|{#1}\right\|}}

\def\R{{\mathbb{R}}}

\newcommand{\ps}{\psi^{(k,c)}_{N,m}}
\newcommand{\jpeven}{P^{(0,k+m/2-1)}_{n}}
\newcommand{\jpodd}{P^{(0,k+m/2)}_{n}}
\newcommand{\wL}{\widetilde{C_{n,m}^0}(Y^\ell_k(x))}
\newcommand{\Le}{C_{n,m}^0(Y^\ell_k(x))}
\journal{Journal of Approximation Theory}

\begin{document}
	
	\begin{frontmatter}
		
		\title{Expansion into Clifford Prolate Spheroidal Wave Functions}
		
		\author[UON]{Hamed Baghal Ghaffari}
		\ead{hamed.baghalghaffari@newcastle.edu.au}
		
		\author[Bizerte]{Ahmed Souabni}
		\ead{souabniahmed@bizerte.r-iset.tn}
		
		\address[UON]{School of Mathematical and Physical Sciences, University of Newcastle, Callaghan NSW 2308, Australia}
		\address[Bizerte]{Faculty of Sciences of Bizerte, University of Carthage, Tunisia}
		
		\begin{abstract}
			In this paper, we investigate the properties of Clifford prolate spheroidal wave functions (CPSWFs) through their associated eigenvalues. We prove that the expansion coefficients in CPSWFs series decay as both the order and the homogeneity degree increase. By establishing a precise connection between the radial CPSWFs and the eigenfunctions of the finite Hankel transform, we derive explicit and non-asymptotic bounds on the corresponding eigenvalues and transfer the spectral decay estimates to the Clifford setting. Consequently, we obtain super-exponential decay rates for the CPSWF expansion coefficients of band-limited Clifford-valued functions. Numerical experiments illustrate both the accuracy and the efficiency of these approximations.
		\end{abstract}
		
		
	\end{frontmatter}
	
	
	\section{Introduction} 
	Time-limited and band-limited functions are fundamental tools in signal processing. According to Heisenberg’s uncertainty principle, a signal cannot be time-limited and band-limited simultaneously. A natural question is that among the band-limited signals, find those signals that are most concentrated on a given interval. This problem was first studied in depth by Slepian, Pollak, and Landau in their pioneering works of the 1960s, where prolate spheroidal wave functions (PSWFs) were introduced as its solutions \cite{slepian1961prolate, landau1961prolate, landau1962prolate, slepian1964prolate}. \\
    Prolate spheroidal wave functions occupy a central place in approximation theory due to their extreme concentration properties and their role as optimally band-limited bases. \\
    Our focus is on the extension proposed by Slepian in \cite{slepian1964prolate}, where the time– frequency concentration problem was formulated in an $m$-dimensional setting. Despite its theoretical importance, the study of signal concentration over bounded domains in higher dimensions has remained relatively underexplored. Nonetheless, several recent contributions address aspects of this topic. In \cite{Lederman}, numerical approaches were proposed for computing multidimensional prolate spheroidal wave functions. Likewise, \cite{Greengard} introduces algorithms for their evaluation, including quadrature rules adapted to band-limited functions and numerical schemes for computing the eigenvalues of the multidimensional finite Fourier transform. 
    Further developments related to multidimensional prolate spheroidal wave functions in the setting of Clifford algebra can be found in \cite{ghaffari2019clifford, propertiesofcliffordlegendre, ghaffari2022higher}. Although the underlying algebraic framework differs substantially from the classical context, the construction of Clifford PSWFs follows a methodology similar to that employed in \cite{Wang3} and \cite{slepian1964prolate}. In particular, \cite{ghaffari2019clifford} establishes a Bonnet-type recurrence relation for Clifford Legendre polynomials, which constitutes a key ingredient in the construction presented in \cite{propertiesofcliffordlegendre}. Additional properties and results are discussed in \cite{ghaffari2022higher}.
	 This extension is natural because classical PSWFs are mainly designed for scalar and one-dimensional settings, and therefore do not fully reflect the geometric structure of higher-dimensional data. Clifford analysis provides a convenient way to handle multicomponent and vector-valued functions. In this context, studying PSWFs allows one to build concentration operators that better respect the geometry of the problem and to obtain basis functions adapted to multidimensional bandlimited signals. For example, the coloured images as higher-dimensional signals using the RGB model can be presented by quaternions (see e.g. \cite{huang2023review, dizon2024holistic} ), as a result, the Clifford PSWFs, defined on dimension $2$, will be able to reconstruct coloured images. 
	
	The present work builds upon the results of \cite{ghaffari2019clifford, ghaffari2022higher, baghal2021clifford, souabni24}, focusing on the analysis of the approximation quality of band-limited functions in the context of Clifford analysis. In this work, we develop a detailed approximation theory for the Clifford prolate spheroidal wave functions (CPSWFs), extending several classical results on multidimensional PSWFs to the Clifford setting. \\
    \noindent{\bf Main contributions. } 
    Our first contribution is a precise spectral analysis of the CPSWFs through an explicit connection with the eigenfunctions of the finite Hankel transform. This link allows us to transfer the known sharp bounds on the Hankel eigenvalues to the Clifford framework to obtain new, non-asymptotic estimates for the eigenvalues of the finite Clifford Fourier transform. Building on these estimates, we establish super-exponential decay rate for the CPSWFs expansion coefficients of band-limited Clifford-valued functions. As a consequence, we derive explicit error bounds for truncated CPSWFs expansions and prove that the resulting approximations converge significantly faster than the classical Fourier–Bessel series. Finally, numerical experiments illustrate the sharpness of our theoretical results and highlight the efficiency of the CPSWFs basis for approximating multidimensional band-limited Clifford-valued signals. \\
    \medskip
    \noindent {\bf Organization of the paper. }
	In Section 2, we review fundamental concepts from Clifford analysis relevant to the developments presented in this work. In Section 3, we give further properties of CPSWFs. The main result here is that we can link the radial CPSWFs and their corresponding eigenvalues to a well-known case in the literature. In Section 4, we investigate the approximation accuracy and speed of convergence of band-limited functions using a series of expansions in terms of CPSWFs, following the approach developed in \cite{souabni24}. The last section presents numerical experiments designed to illustrate the theoretical findings of this study. These simulations provide an assessment of the accuracy of the proposed approximations.

	\section{Background}
	Let $\mathbb{R}^{m}$ be $m$-dimensional Euclidean space and let
	$\{e_{1},e_{2},\dots ,e_{m}\}$
	be an orthonormal basis for
	$\mathbb{R}^{m}.$
	We endow these vectors with the multiplicative properties
	\begin{align*}
		e_{j}^{2}&=-1,\; \; j=1,\dots , m,\\
		e_{j}e_{i}&=-e_{i}e_{j}, \;\; i\neq j, \;\; i,j=1,\dots , m.
	\end{align*}
	For any subset
	$A=\{j_{1},j_{2},\dots, j_{h}\}\subseteq \{1,\dots ,	m\}=Q_m,$ with $j_1<j_2<\cdots <j_h$
	we consider the formal product
	$e_{A}=e_{j_{1}}e_{j_{2}}\dots e_{j_{h}}.$
	Moreover, for the empty set
	$\emptyset$
	one puts
	$e_{\emptyset}=1$ (the identity element). The Clifford algebra ${\mathbb R}_m$ is then the $2^m$-dimensional real associative algebra 
	$${\mathbb R}_m=\bigg\{\sum\limits_{A\subset Q_m}\lambda_Ae_A:\, \lambda_A\in{\mathbb R}\bigg\}.$$
	Similarly, the Clifford algebra ${\mathbb C}_m$ is the $2^m$-dimensional complex associative algebra
	$${\mathbb C}_m=\bigg\{\sum\limits_{A\subset Q_m}\lambda_Ae_A:\, \lambda_A\in{\mathbb C}\bigg\}.$$
	Every element $\lambda=\sum\limits_{A\subset Q_m}\lambda_Ae_A\in{\mathbb C}_m$ may be decomposed as  
	$\lambda=\sum\limits_{k=0}^{m}[\lambda]_{k},$
	where 
	$[\lambda]_{k}=\sum\limits_{\vert A\vert=k}\lambda_{A}e_{A}$
	is the so-called 
	$k$-vector
	part of 
	$\lambda\, (k=0,1,\dots ,m).$
	
	Denoting by 
	$\mathbb{R}_{m}^{k}$
	the subspace of all 
	$k$-vectors
	in
	$\mathbb{R}_{m},$
	i.e., the image of 
	$\mathbb{R}_{m}$
	under the projection operator 
	$[\cdot]_{k},$
	one has the multi-vector decomposition
	$\mathbb{R}_{m}=\mathbb{R}_{m}^{0}\oplus \mathbb{R}_{m}^{1}\oplus\cdots \oplus \mathbb{R}_{m}^{m},$
	leading  to the identification of
	$\mathbb{R}$
	with the subspace of real scalars
	$\mathbb{R}_{m}^{0}$
	and of
	$\mathbb{R}^{m}$
	with the subspace of real Clifford vectors 
	$\mathbb{R}_{m}^{1}.$ The latter identification is achieved by identifying the point
	$(x_{1},\dots,x_{m})\in{\mathbb R}^m$
	with the Clifford number
	$x=\sum\limits_{j=1}^{m}e_{j}x_{j}\in{\mathbb R}_m^1$.
	The Clifford number 
	$e_{M}=e_{1}e_{2}\cdots e_{m}$
	is called the pseudoscalar; depending on the dimension 
	$m,$
	the pseudoscalar commutes or anti-commutes with the 
	$k$-vectors
	and squares to 
	$\pm 1.$
	The Clifford conjugation on ${\mathbb C}^m$ is the conjugate linear mapping $\lambda\mapsto\bar\lambda$ of ${\mathbb C}_m$ to itself satisfying
	\begin{align*}
		\overline{\lambda \mu}&=\bar\mu\bar{\lambda},\;\;\;\; \textnormal{for all}\;\lambda,\mu\in\mathbb{C}_{m},\\
		\overline{\lambda_{A}e_{A}}&=\overline{\lambda_{A}}\overline{e_{A}},\;\;\; \lambda_A\in\mathbb{C},\\
		\overline{e_{j}}&=-e_{j},\;\; j, \;\; j=1,\cdots , m.
	\end{align*}
	The Clifford conjugation leads to a Clifford inner product $\langle \cdot,\cdot\rangle$ and an associated norm $|\cdot |$ on 
	$\mathbb{C}_{m}$
	given respectively by
	$$\langle \lambda,\mu\rangle=[\bar{\lambda}\mu]_{0}\;\;\;\textnormal{and}\;\;\; \vert\lambda\vert^{2}=[\bar{\lambda}\lambda]_{0}=\sum\limits_{A}\vert\lambda_{A}\vert^{2},$$
	for $\lambda=\sum_{A\subset Q_M}\lambda_Ae_A\in{\mathbb C}_m$. 
	
	The product of two vectors $x$, $y\in{\mathbb R}_m^1$ can be decomposed as the sum of a scalar and a 2-vector, also called a bivector:
	$$xy=-\langle x, y\rangle +x\wedge y,$$
	where
	$\langle x,y\rangle=-\sum\limits_{j=1}^{m}x_{j}y_{j}\in \mathbb{R}^{0}_{m}$,
	and,
	$x\wedge y=\sum\limits_{i=1}^{m}\sum\limits_{j=i+1}^{m}e_{i}e_{j}(x_{j}y_{j}-x_{j}y_{i})\in\mathbb{R}^{2}_{m}$.
	Note that the square of a vector 
	$x$
	is scalar-valued and equals the norm squared up to the minus sign:
	$$x^{2}=-\langle x,x\rangle=-\vert x\vert^{2}.$$
	Clifford analysis offers a function theory that is a higher-dimensional analogue of the theory of holomorphic functions of one complex variable. The functions considered are defined in the Euclidean space 
	$\mathbb{R}^{m}$
	and take their values in the Clifford algebra 
	$\mathbb{R}_{m}.$
	
	The central notion in Clifford analysis is monogenicity, a multidimensional counterpart to holomorphy in the complex plane.
	\begin{defn}
		A function 
		$f(x)=f(x_{1},\dots, x_{m})$
		defined and continuously
		differentiable in an open region of 
		$\mathbb{R}^{m}$
		and taking values in 
		$\mathbb{C}_{m}$
		is said to be left monogenic in that region if 
		$$\partial_{x}f=0,$$
		where
		$\partial_{x}=\sum_{j=1}^me_j\partial_{x_j}$
		is the Dirac operator 
		and 
		$\partial_{x_{j}}$
		is the partial differential operator
		$\frac{\partial}{\partial x_{j}}.$
		We also define the Euler differential operator by
		$E=\sum_{j=1}^{m}x_{j}\partial_{x_{j}}$.
		The Laplace operator is factorized by the Dirac operator as follows:
		\begin{equation}
			\Delta_{m}=-\partial_{x}^{2}.
		\end{equation}
	\end{defn}
	The notion of right monogenicity is defined similarly by letting the Dirac operator act from the right. A ${\mathbb C}_m$-valued function $f(x)=\sum_{A\subset Q_m}f_A(x)e_A$ (where each $f_A$ takes complex values) is left monogenic if and only if its Clifford conjugate $\bar{f}(x)=\sum_{A\subset Q_m}\overline{f_A}(a)\overline{e_A}$
	is right monogenic. In fact, $\overline{\partial f}=-\overline{f}\partial$.
	\begin{defn}\label{monogenic_homogeneous_polynomial_Def}
		A left (respectively right) monogenic homogeneous polynomial
		$P_{k}$
		of degree 
		$k\; (k\geq 0)$
		in
		$\mathbb{R}^{m}$
		is called a left (respectively right) solid inner spherical monogenic of order 
		$k.$
		The set of all left (respectively right) solid inner spherical monogenic of order 
		$k$
		will be denoted by
		$M_{l}^{+},$
		respectively
		$M_{r}^{+}.$
		It can be shown 
		\cite{delanghe2012clifford}
		that the dimension of 
		$M_{l}^{+}(k)$
		is given by 
		$$\dim M_{l}^{+}(k)=\frac{(m+k-2)!}{(m-2)!k!}=d_{k,m}.$$
		A left (respectively right) monogenic homogeneous function $Q_k$ of degree $-(k+m-1)$ in ${\mathbb R}^m \setminus\{0\}$ is called a left (respectively right) solid outer spherical monogenic of order $k.$
	\end{defn}
	\begin{lem}
		We can see that 
		$$M_{l}^{+}(k)\cap M_{l}^{-}(k)=\{0\}.$$
	\end{lem}
	We let $M_{l}(k):=M_{l}^{+}(k)\oplus M_{l}^{-}(k)$
	and note that if $P_{k}\in M_{l}^{+}(k)$ and $Q_{k}(x)=\dfrac{x}{\vert x\vert^{m}}P_{k}(\dfrac{x}{\vert x\vert^2})$ Then $Q_k\in M_{l}^{-}(k)$. 
	
	\begin{lem}\cite{delanghe2012clifford}\label{diracderivativelemma}
		If
		$P_{k}\in M_{l}^{+}(k)$
		and 
		$s\in \mathbb{N}$, then 
		\begin{equation*}
			\partial_{x}[x^{s}P_{k}(x)]=
			\begin{cases}
				-sx^{s-1}P_{k}(x)&\text{ if $s$ is even},\\
				-(s+2k+m-1)x^{s-1}P_{k}(x)&\text{ if $s$ is odd.}
			\end{cases}
		\end{equation*}
	\end{lem}
	\begin{proof}
		For the proof see \cite{delanghe2012clifford}.
	\end{proof}
	\begin{defn}
		A real-valued polynomial
		$S_{k}$
		of degree
		$k$
		on
		$\mathbb{R}^{m}$
		satisfying
		$$\Delta_{m}S_{k}(x)=0,\;\;\;\;\; \textnormal{and,}\;\;\;\;\; S_{k}(tx)=t^{k}S_{k}(x)\quad (t>0),$$
		is called a solid spherical harmonic of degree $k$. The collection of solid spherical harmonics of degree $k$ on
		$\mathbb{R}^{m}$
		is denoted
		$\mathcal{H}(k)$
		(or $\mathcal{H}(m,k)$).
	\end{defn}
	Since
	$\partial^{2}_{x}=-\Delta_{m}$,
	we have that
	$$M_{l}^{+}(k)\subset \mathcal{H}(k),\;\;\;\;\; \textnormal{and,}\;\;\;\;\; M_{r}^{+}(k)\subset \mathcal{H}(k).$$
	Let
	$H_{(r)}$
	be a unitary right Clifford-module, i.e. 
	$(H_{(r)},+)$
	is an abelian group
	and a law of scalar multiplication
	$(f,\lambda)\to f\lambda$
	from
	$H_{(r)}\times\mathbb{C}_{m}$
	into
	$H_{r}$
	is defined such that for all
	$\lambda,\mu\in\mathbb{C}_{m}$
	and
	$f,g\in H_{(r)}:$
	\begin{align*}
		&(i)\;f(\lambda+\mu)=f\lambda+f\mu,\hspace*{10cm}\\
		&(ii)\;f(\lambda\mu)=(f\lambda)\mu,\\
		&(iii)\;(f+g)\lambda=f\lambda+g\lambda,\\
		&(iv)\;fe_{\emptyset}=f.
	\end{align*}
	Note that
	$H_{(r)}$
	becomes a complex vector space if
	$\mathbb{C}$
	is identified with
	$\mathbb{C}e_{\emptyset}\subset\mathbb{C}_{m}.$
	Then a function
	$\langle\cdot,\cdot\rangle:H_{(r)}\times H_{(r)}\to\mathbb{C}_{m}$
	is said to be a ${\mathbb C}_m$-valued  inner product on
	$H_{(r)}$
	if for all
	$f,g,h\in H_{(r)}$
	and
	$\lambda\in\mathbb{C}_{m}:$
	\begin{align*}
		&(i)\;\langle f,g\lambda+h\rangle=\langle f,g\rangle\lambda+\langle f,h\rangle,\hspace*{10cm}\\
		&(ii)\;\langle f,g\rangle=\overline{\langle g,f\rangle},\\
		&(iii)\;[\langle f,f\rangle]_{0}\geq 0 \;\; \textnormal{and}\;\; [\langle f,f\rangle]_{0}= 0 \;\;\textnormal{if and only if}\;\; f=0.
	\end{align*}
	From this
	${\mathbb C}_m$-valued
	inner product
	$(\cdot,\cdot)$,
	one can recover the complex inner product
	$$( f,g) =[\langle f,g\rangle]_{0},$$
	on
	$H_{r}$. Putting for each
	$f\in H_{(r)}$,
	\begin{equation}
		\Vert f\Vert^{2}=( f,f) ,\label{ip norm}
	\end{equation}
	$\Vert \cdot\Vert$
	becomes a norm on 
	$H_{r}$
	turning it into a normed right Clifford-module. 
	
	Let
	$H_{(r)}$
	be a unitary right Clifford-module provided with a ${\mathbb C}_m$-valued inner product
	$\langle\cdot,\cdot\rangle.$
	Then it is called a right Hilbert Clifford-module if
	$H_{(r)}$
	considered as a complex vector space provided with the complex inner product 
	$(\cdot,\cdot)$
	is a Hilbert space. 
	
	We consider the ${\mathbb C}_m$-valued inner product of the functions 
	$f,g:{\mathbb R}^m\to{\mathbb C}_m$ defined by
	by
	$ \langle f,g\rangle=\int\limits_{\mathbb{R}^{m}}\overline{f(x)}g(x)\, dx$,
	where
	$dx$
	is the Lebesgue measure on 
	$\mathbb{R}^{m}$. 
	The associated norm is given by (\ref{ip norm}).
	The right Clifford-module of ${\mathbb C}_m$-valued measurable functions on
	$\mathbb{R}^{m}$
	for which 
	$\Vert f\Vert^{2}<\infty$ is denoted 
	$L^{2}(\mathbb{R}^{m},{\mathbb C}_m)$.
	
	The standard tensorial multi-dimensional Fourier transform is given by
	\begin{equation*}
		\mathcal{F}f(\xi)=\int\limits_{\mathbb{R}^{m}}e^{-2\pi i\langle x,\xi\rangle}f(x)\, dx
	\end{equation*}
	for $f\in L^1({\mathbb R}^m,{\mathbb C}_m)$.
	\begin{defn}
		For any $f,g\in L^{2}(\mathbb{R}^{m},{\mathbb C}_m),$ the operator linear operator $T: L^2({\mathbb R}^m,{\mathbb C}_m)\to L^2({\mathbb R}^m,{\mathbb C}_m)$ is self-adjoint if 
		$$\langle Tf,g\rangle=\langle f,Tg\rangle.$$
	\end{defn}
	
	\begin{defn}
		The linear operator $T$ is compact, if $\{T(f_{n})\}$ has a convergent subsequence for every bounded sequence of $\{f_{n}\}\in L^{2}(\mathbb{R}^{m},{\mathbb C}_m)$.
	\end{defn}
	
	\begin{thm}\label{planchereltheorem}
		The Fourier transform 
		$\mathcal{F}$
		extends to  an isometry on $L^2({\mathbb R}^m,{\mathbb C}_m)$, i.e., for all
		$f,g\in L^{2}(\mathbb{R}^{m},{\mathbb C}_m)$ the Parseval formula holds:
		$$\langle f,g\rangle=\langle \mathcal{F}f,\mathcal{F}g\rangle.$$
		In particular, for each 
		$f\in L_{2}(\mathbb{R}^{m},{\mathbb C}_m)$
		one has $\Vert f\Vert_{2}=\Vert \mathcal{F}f\Vert_{2}$.
	\end{thm}
	\begin{prop}
		For any two monogenic, homogeneous polynomials 
		$Y_{k}$ of homogeneous degree $k$
		and 
		$Y_{l}$ of homogeneous degree $l$
		we have 
		$$\langle Y_{k},Y_{l}\rangle_{L^{2}(S^{m-1})}:=\int\limits_{S^{m-1}}\overline{Y_k(\omega )}Y_l (\omega )\, d\omega=0$$
		if $k\neq l$.
	\end{prop}
	
	The following result may be obtained as a simple application of the Clifford-Stokes theorem.
	\begin{lem}\label{property of two monogenic and x between}
		Let $f$, $g$ be defined on a neighbourhood $\Omega $ of the unit ball in ${\mathbb R}^m$ and suppose $f$ is right monogenic on $\Omega$ while $g$ is left monogenic on $\Omega$.
		Then
		\begin{equation}
			\int_{S^{m-1}}f(\omega )\omega g(\omega )\, d\omega =0.
		\end{equation}
	\end{lem}
	The following well-known result appears as Lemma 9.10.2 in \cite{andrews1999special} and is a corollary of the Funk-Hecke Theorem \cite{groemer1996geometric}.
	\begin{lem}\label{Lemma_from_Chinese_Paper}	
		Let
		$\hat{\xi},\omega\in S^{m-1}$, $r>0$ and $Y_k\in \mathcal{H}_{k}^{m}$. Then 
		$$\int_{S^{m-1}}e^{-2\pi ir\langle \hat{\xi},\omega \rangle}Y_{k}(\omega )\,d\sigma(\omega )=\frac{2\pi(-i)^{k}}{r^{\frac{m}{2}-1}}J_{k+\frac{m}{2}-1}(2\pi r)Y_{k}(\hat{\xi})$$
		where
		$J_{k+\frac{m}{2}-1}$
		is a Bessel function of the first kind.
	\end{lem}
	
	As a consequence of Lemma \ref{Lemma_from_Chinese_Paper}, we have for all $Y_k\in \mathcal{H}_{k}^{m}$
	\begin{equation}
		\int_{B(1)}e^{-2\pi i\langle x,\xi\rangle}Y_k(x)\, dx=\frac{(-i)^k}{|\xi |^{\frac{m}{2}+k}}J_{\frac{m}{2}+k}(2\pi |\xi |)Y_k(\xi ).\label{solid F-H thm}
	\end{equation}
	For $r>0$, let $B(r)$ be the closed ball of radius $r$ and centre $0$ in ${\mathbb R}^m$, i.e.,
	$$B(r)=\{x\in{\mathbb R}^m:\, |x|\leq r\}.$$
	\begin{defn}\label{Definition of Gegenbauer}
		The Clifford-Legendre polynomial $C_{n,m}^{0}(Y_{k})(x)$ can be defined as 
		$$\Le=\partial_{x}^{n}((1-\vert x\vert^{2})^{n}Y_{k}(x)).$$
	\end{defn}
	We denote by $\wL$ the normalized Clifford Legendre polynomial so that their 2-norm is equal to 1: $$ \wL(x) = \frac{1}{h_{k,n}} \Le(x) \quad h_{k,n} = \frac{\sqrt{2k+2n+m}}{2^n n!},$$
	where $n$ is the order, $k$ is the homogeneity degree and $m$ is the dimension. 
	
	The Clifford Legendre polynomial can be written in the following form:
	$$
	\begin{array}{rcl}
		C^0_{2n,m}(Y^\ell_k(x)) &=& (-1)^n2^{2n} (2n)! \, \jpeven(2|x|^2-1) Y^{(\ell)}_k(x) \\
		C^0_{2n+1,m}(Y^\ell_k(x)) &=& (-1)^{n+1}2^{2n} (2n+1)!\, x\, \jpodd(2|x|^2-1)  Y^{(\ell)}_k(x),
	\end{array}
	$$
	where $P^{(\alpha,\beta)}_n$ are the Jacobi polynomials  
	$$ P_{n}^{(\alpha ,\beta)}(x) =\frac{(-1)^{n}}{2^{n}n!}(1-x)^{-\alpha
	}(1+x)^{-\beta }\frac{d^{n}}{dx^{n}}\left[(1-x)^{n+\alpha }(1+x)^{n+\beta
	}\right]. $$
	
	\begin{defn}
		Given $c\geq 0$, the Clifford differential operator 
		$L_{c}$
		acts of 
		$C^{2}(B(1),\mathbb{R}_{m})$
		as follows:
		\begin{equation}\label{Definition_of_Lc_operator}
			L_{c}f(x)=\partial_{x}((1-\vert x\vert^{2})\partial_{x}f(x))+4\pi^{2}c^{2}\vert x\vert^{2}f(x)
		\end{equation}
		where 
		$B(1)$
		is the unit ball in
		$\mathbb{R}^{m}$
		and 
		$\partial_{x}$
		is the Dirac Operator. We define the Clifford Prolate Spheroidal Wave Functions (CPSWFs) as the eigenfunctions of 
		$L_{c}.$
	\end{defn}
	\begin{prop}\label{Lc_Self_Adjoint}
		The operator $L_{c}$ defined in 
		\eqref{Definition_of_Lc_operator}
		is self-adjoint.
	\end{prop}
	\noindent \textbf{Consequence :} 
	Given $c>0$, the CPSWFs are equivalently the eigenfunctions of the finite Fourier transformation $\mathcal{G}_c$. So we have 
	$$ \int_{\overline{B(1)}}\psi_{N,m}^{k,c}(y)e^{2\pi i c \langle x,y \rangle}\, dy=\mu_{N,m}^{k,c}\psi_{N,m}^{k,c}(x).$$
	

	\section{Estimation of the eigenvalues}
	This section will be devoted to proving several key results that are essential for our proof of the development of a function in the prolate basis.
	First, we will provide an upper and lower bound for the eigenvalues corresponding to the differential operator $L_c$.
	\begin{prop}\label{Chi_Relation_Theorem}
		Let $n$ be an even number and the real constants $\chi_{n,m}^{k,0}$, $\chi_{n,m}^{k,c}$, be given by
		\begin{align*}
			\wL&=\chi_{n,m}^{k,0}\, \wL,\\
			L_c\psi_{n,m}^{k,c,i}&=\chi_{n,m}^{k,c}\psi_{n,m}^{k,c,i}.
		\end{align*}
		Then for $c>0$ we have that 
		\begin{equation}\label{Chi_Relation}
			n(n+2k+m)<\chi_{n,m}^{k,c}<n(n+2k+m)+8\pi^2 c^2
		\end{equation}
	\end{prop}
	\begin{proof}
		By taking derivatives in terms of $c$ from \eqref{Definition_of_Lc_operator} we have that 
		$$
		L_{c} \partial_{c}\psi_{n,m}^{k,c}(x)	+ 4\pi^2 2c\vert x\vert^{2} \psi_{n,m}^{k,c}=\chi_{n,m}^{k,c}  \partial_{c}\psi_{n,m}^{k,c}(x) +\partial_{c} \chi_{n,m}^{k,c}  \psi_{n,m}^{k,c}(x),
		$$
		hence 
		$$ \langle (L_{c}-\chi_{n,m}^{k,c}) \partial_{c}\psi_{n,m}^{k,c}(x) ,  \psi_{n,m}^{k,c}(x) \rangle= \langle (\partial_{c} \chi_{n,m}^{k,c} - 8\pi^2 c\vert x\vert^{2}  ) \psi_{n,m}^{k,c}(x) ,  \psi_{n,m}^{k,c}(x) \rangle. $$
		By Proposition \ref{Lc_Self_Adjoint} we have 
		$$ \langle  \partial_{c}\psi_{n,m}^{k,c}(x) ,  (L_{c}-\chi_{n,m}^{k,c}) \psi_{n,m}^{k,c}(x) \rangle= \langle (\partial_{c} \chi_{n,m}^{k,c} - 8\pi^2 c\vert x\vert^{2}  ) \psi_{n,m}^{k,c}(x) ,  \psi_{n,m}^{k,c}(x) \rangle $$
		hence, 
		{\small
			\begin{equation}\label{Auxilary_Chi_Proof1}
				0=\langle (\partial_{c} \chi_{n,m}^{k,c} - 8\pi^2 c\vert x\vert^{2}  ) \psi_{n,m}^{k,c}(x) ,  \psi_{n,m}^{k,c}(x) \rangle=\partial_{c}\chi_{n,m}^{k,c}-8\pi^{2}c\int_{B(1)}\vert x\vert^{2}  \psi_{n,m}^{k,c}(x)  \psi_{n,m}^{k,c}(x) dx.
		\end{equation}}
		By \eqref{Auxilary_Chi_Proof1} we have
		\begin{equation}\label{Auxilary_Chi_Proof2}
			0\leq \partial_{c}\chi_{n,m}^{k,c}-8\pi^2 c^2 \int_{B(1)}\vert x\vert^{2}  \psi_{n,m}^{k,c}(x)  \psi_{n,m}^{k,c}(x) dx\leq 8\pi^2 c.
		\end{equation}
		Now by \eqref{Auxilary_Chi_Proof2} we get 
		$$ \frac{\chi_{n,m}^{k,c}-\chi_{n,m}^{k,0}}{c}< 8\pi^{2}c $$
		which concludes the\eqref{Chi_Relation}.
	\end{proof}
	The following result establishes a relationship between the Clifford PSWFs and the Hankel prolate functions. The advantage of this result lies in the fact that the eigenvalues associated with this particular class of prolate functions (Hankel or circular) have been extensively studied in the literature (see, e.g., \cite{non_asympt_mourad,BK2017,Karoui_Moumni}). First, let us establish some concepts. 
	\begin{defn}
		For a positive real number $c$, we denote by $\mathcal{H}^{(\alpha)}_c $ the finite Hankel transform given by 
		$$ \mathcal{H}^{(\alpha)}_c.f(x) =\sqrt{2\pi} \int_{0}^1 \sqrt{2\pi c xy} J_\alpha( 2\pi c xy) f(y)dy.$$
		We denote $\varphi^{(\alpha)}_n$ the eigenfunction of $\mathcal{H}^{(\alpha)}_c $ and $\gamma^{(\alpha)}_n$ the associated eigenvalue.
	\end{defn}
	
	\begin{prop}\label{relation_circular_case}
		Let $c>0$ and $k,m\in\mathbb{N}$.  
		Denote the even and odd Clifford prolate spheroidal wave functions (CPSWFs) by
		\[
		\psi^{(k,c)}_{2N,m}(x)
		= P^{k,c}_{N,m}(|x|^{2})\, Y_k(x), 
		\qquad 
		\psi^{(k,c)}_{2N+1,m}(x)
		= Q^{k,c}_{N,m}(|x|^{2})\, x\, Y_k(x),
		\]
		where $Y_k$ is a spherical monogenic of degree $k$.

		Then:
		
		\begin{itemize}
			\item For even CPSWFs,
			\[
			\gamma^{\left(k+\frac{m}{2}-1\right)}_N
			= (-i)^k\, c^{\frac{m-1}{2}}\mu^{(k,c)}_{2N,m},
			\qquad
			\varphi^{\left(k+\frac{m}{2}-1\right)}_N(r)
			= r^{\,k+\frac{m-1}{2}}\, P^{k,c}_{N,m}(r^2).
			\]
			
			\item For odd CPSWFs,
			\[
			\gamma^{\left(k+\frac{m}{2}\right)}_N
			= (-i)^{k+1}\, c^{\frac{m-1}{2}} \mu^{(k,c)}_{2N+1,m},
			\qquad
			\varphi^{\left(k+\frac{m}{2}\right)}_N(r)
			= r^{\,k+\frac{m+1}{2}}\, Q^{k,c}_{N,m}(r^2).
			\]
		\end{itemize}
	\end{prop}
	\begin{proof}
		We know that the CPSWFs are the eigenfunctions of the finite Fourier transformations. So we have that
		$$\mathcal{G}_{c}\psi_{2N,m}^{k,c}(x)=\mu_{2N,m}^{k,c}\psi_{2N,m}^{k,c}(x),$$
		which means
		$$ \int_{\overline{B(1)}}\psi_{2N,m}^{k,c}(y)e^{2\pi i c \langle x,y \rangle}\, dy=\mu_{2N,m}^{k,c}\psi_{2N,m}^{k,c}(x),$$
		therefore, we have that
		\begin{equation}\label{even case integral opertor in polar}
			\int_{S^{m-1}}\int_{0}^{1}P_{N}^{k}(r^2)r^{k+m-1}Y_{k}(\omega) e^{2\pi i c s r\langle \omega,\xi \rangle}\, dr d\omega=\mu_{2N}^{k,c}P_{N,m}^{k}(s^2)s^{k}Y_{k}(\omega),
		\end{equation}
		But we may change the order of the integration
		{ \tiny
			\begin{equation} \label{even case integral opertor in polar change int order}
				\int_{S^{m-1}}\int_{0}^{1}P_{N}^{k}(\vert x\vert^2)r^{k+m-1}Y_{k}(\omega) e^{2\pi i c s r\langle \omega,\xi \rangle}\, dr d\omega=\int_{0}^{1}(\int_{S^{m-1}}Y_{k}(\omega)e^{2\pi i c s r\langle \omega,\xi \rangle}d\omega)P_{N}^{k}(r^2)r^{k+m-1} dr.
			\end{equation} 
		}
		Now, we need to calculate the inner integral. From Lemma 9.10.2 in \cite{andrews1999special}. We know that 
		$$\int_{S^{m-1}}Y_{k}(\omega)e^{-2\pi i r\langle \omega,\hat{\xi} \rangle}d\omega=\frac{2\pi (-i)^{k}}{r^{\frac{m}{2}-1}}J_{k+\frac{m}{2}-1}(2\pi r)Y_{k}(\hat{\xi}).$$
		So we have that 
		\begin{align}\label{the formula obtained for even case}
			\int_{S^{m-1}}Y_{k}(\omega)e^{2\pi i c s r\langle \omega,\xi \rangle}d\omega=&\int_{S^{m-1}}Y_{k}(-\omega)e^{-2\pi i (c s r)\langle \omega,\xi \rangle}d\omega\hspace*{7cm}\nonumber\\
			=&(-1)^{k}\int_{S^{m-1}}Y_{k}(\omega)e^{-2\pi i (c s r)\langle \omega,\xi \rangle}d\omega\nonumber\\
			=&\frac{(-1)^{k}2\pi(-i)^{k}}{(csr)^{\frac{m}{2}-1}}J_{k+\frac{m}{2}-1}(2\pi c s r)Y_{k}(\xi)\nonumber\\
			=&\frac{2\pi i^{k}}{(csr)^{\frac{m}{2}-1}} J_{k+\frac{m}{2}-1}(2\pi c s r)Y_{k}(\xi).
		\end{align}
		Therefore, the integral operator for the radial part using \eqref{even case integral opertor in polar change int order}, \eqref{even case integral opertor in polar}, and, \eqref{the formula obtained for even case} can be obtained as
		$$\mu_{2N,m}^{k,c}P_{N,m}^{k,c}(s^2)s^k Y_{k}(\xi)=\int_{0}^{1}\frac{2\pi i^{k}}{(csr)^{\frac{m}{2}-1}}J_{k+\frac{m}{2}-1}(2\pi c s r) P_{N,m}^{k,c}(r^{2}) r^{k+m-1} dr\; Y_{k}(\xi),$$
		in which by cancellation $Y_{k}(\xi)$ we get
		$$\mu_{2N,m}^{k,c}P_{N,m}^{k,c}(s^2)s^{k+\frac{m}{2}-1}=\frac{2\pi i^{k}}{c^{\frac{m}{2}-1}}\int_{0}^{1}J_{k+\frac{m}{2}-1}(2\pi c s r) P_{N,m}^{k,c}(r^{2}) r^{k+\frac{m}{2}} dr,$$
		now if we set $S_{N,m,even}^{k,c}(r)=P_{N,m}^{k,c}(r^{2}) r^{k+\frac{m-1}{2}}$, then we have 
		\begin{equation}\label{integral operator for even radial prolate}
			\mu_{2N,m}^{k,c}S_{N,m,even}^{k,c}(s)=\frac{\sqrt{2\pi} i^{k}}{c^{\frac{m-1}{2}}}\int_{0}^{1}\sqrt{2\pi csr}J_{k+\frac{m}{2}-1}(2\pi c s r) S_{N,m,even}^{k,c}(r)dr.
		\end{equation}

		Now, we try to obtain a similar formula for the odd case. The process is almost similar. We need to calculate
		$\int_{S^{m-1}}\omega Y_{k}(\omega)e^{2\pi i c r s\langle \omega,\xi \rangle }d\omega.$ Let $x=s\xi$ so we have that
		$$F(x):=\int_{S^{m-1}} Y_{k}(\omega)e^{2\pi i c r \langle \omega,x \rangle }d\omega.$$
		We apply the equation by the Dirac operator from the left. Therefore, we have
		\begin{align*}
			2\pi i c r&\int_{S^{m-1}} \omega Y_{k}(\omega)e^{2\pi i c r \langle \omega,x \rangle }d\omega\\ &=\partial_{x}[\frac{2\pi i^{k}}{(csr)^{\frac{m}{2}-1}}J_{k+\frac{m}{2}-1}(2\pi c s r) Y_{k}(\xi)]\hspace*{4cm}\\
			&=\frac{2\pi i^{k}}{(c r)^{\frac{m}{2}-1}}\sum_{j}e_{j}\frac{\partial}{\partial x_{j}}[J_{k+\frac{m}{2}-1}(2\pi c r\vert x\vert)\vert x\vert^{-(k+\frac{m}{2}-1)}Y_{k}(x)]\\
			&=\frac{2\pi i^{k}}{(c r)^{\frac{m}{2}-1}}\sum_{j}e_{j}\frac{\partial}{\partial x_{j}}[J_{k+\frac{m}{2}-1}(2\pi c r\vert x\vert)\vert x\vert^{-(k+\frac{m}{2}-1)}]Y_{k}(x)\\
			&=\frac{2\pi i^{k}}{(c r)^{\frac{m}{2}-1}}(-2\pi c r)x [J_{k+\frac{m}{2}}(2\pi c r\vert x\vert)\vert x\vert^{-(k+\frac{m}{2})}]Y_{k}(x).
		\end{align*}
		Therefore,
		$$\int_{S^{m-1}} \omega Y_{k}(\omega)e^{2\pi i c r \langle \omega,x \rangle }d\omega=\frac{-(2\pi) i^{k-1}}{(c r)^{\frac{m}{2}-1}}J_{k+\frac{m}{2}}(2\pi c s r) s^{-(\frac{m}{2}-1)} \xi Y_{k}(\xi).$$
		By assuming that $S_{N,m,odd}^{k,c}(r)=Q_{N,m}^{k,c}(r^{2}) r^{k+\frac{m+1}{2}},$ we get a similar integral operator for odd versions
		\begin{equation}\label{integral operator for odd radial prolate}
			\mu_{2N+1,m}^{k,c}S_{N,m,odd}^{k,c}(s)=\frac{\sqrt{2\pi} i^{k+1}}{c^{\frac{m-1}{2}}}\int_{0}^{1} \sqrt{2\pi csr} J_{k+\frac{m}{2}}(2\pi c s r) S_{N,m,odd}^{k,c}(r) dr.
		\end{equation}
	\end{proof}
	\begin{rem}
		In the current proof, in fact, the odd version is the same as the even version with $k+1$. The reason is that for any spherical monogenic, $Y_{k}(x)$, the $x Y_{k}(x)$ is a spherical harmonic of degree $k+1$.
	\end{rem}
	
	\noindent \textbf{Consequence :} Using the result of this last proposition together with Theorems 3.1 and 3.2, in \cite{non_asympt_mourad}, one gets the following estimates of the eigenvalues associated with the Clifford finite Fourier transform on $B(0, c)$, 
	\begin{itemize}
		\item If $N < \frac{c}{2}$, then $ \vert \mu_{N,m}^{k,c}\vert^{2} \geq \frac{1}{c^{m}}\left( 1-\frac{10 \, c^{k+\frac{m}{2}+2N-1}}{N!\, e^{c}}\right)$. 
		\item If $N> \frac{ec}{4}$, then $  \vert \mu_{N,m}^{k,c}\vert^{2} \leq \frac{1}{c^{m}}\left( \frac{e c}{4N+2k+m+3}\right)^{2N+k+\frac{m}{2}} $
	\end{itemize}
	
	\begin{rem}
		The last result is illustrated in the following figure, which shows that the behavior of the eigenvalues consists of three parts. The first shows that $\mu_N$ is close to 1. The second is the plunge region (this region is not studied in this paper). Then, the eigenvalues exhibit super-exponential decay beyond $\frac{ec}{4}$. Visibly, this value is not optimal.
	\end{rem}
	In the following figures, we represent the graph of the eigenvalues for different values of $c$ and $n$.
	\begin{figure}[H]
		\centering
		\includegraphics[width=.7\linewidth]{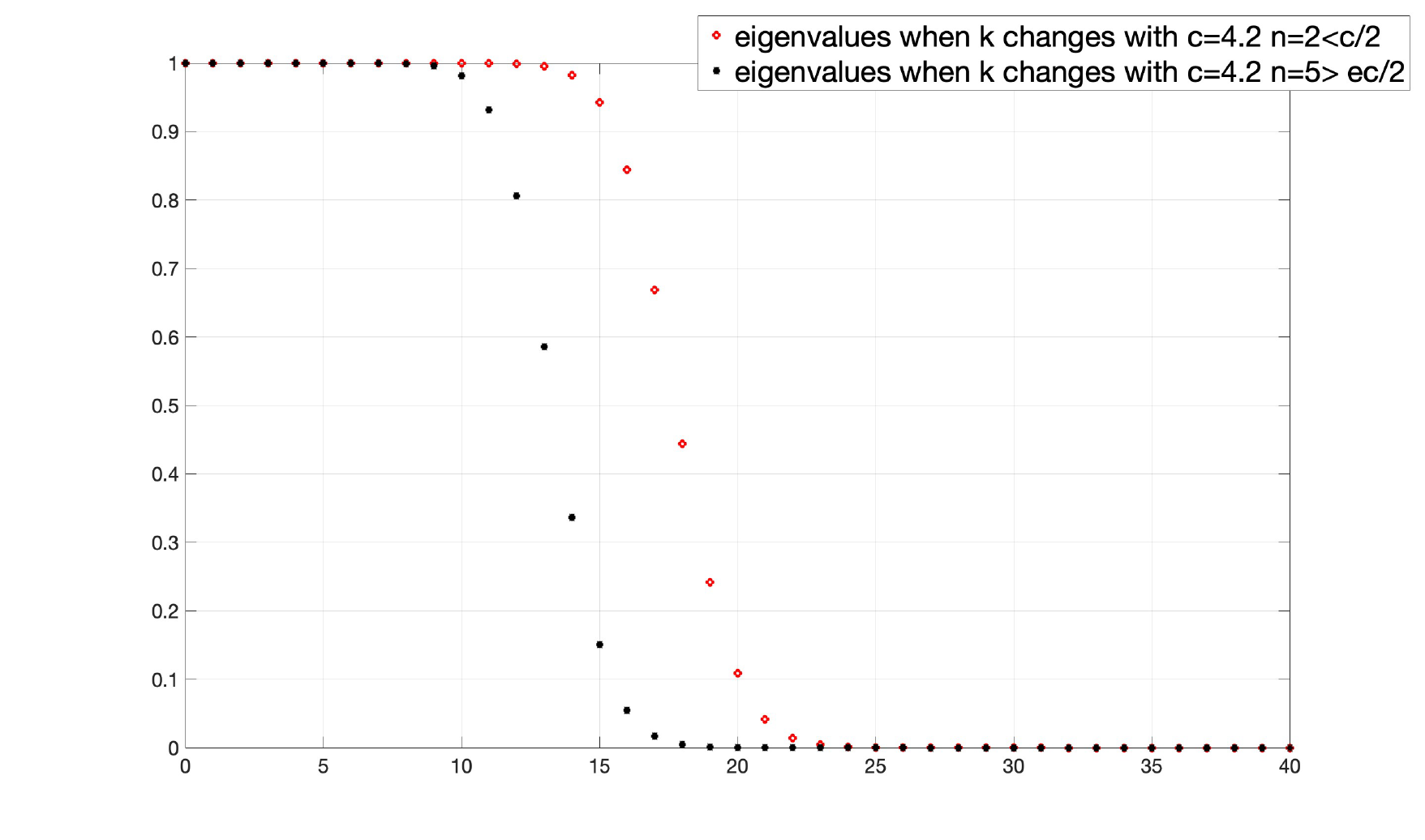}
		\caption{The Graph of Changes of the eigenvalues of the 2-dimension CPSWFs as $k$ changes ($k=0,\cdots, 40$) for different values of the $c=4.2$ and $n=2,5$}\label{equality_of_the_eigenvalues}
	\end{figure}
	
	\section{Expansion into Clifford PSWFs}
	This section aims to study the quality of the approximation in the framework of the d-dimensional Clifford prolate spheroidal wave function series expansion. 
	We denote by $\hat{f}$, the Fourier transform of $f$, 
	
	$$  \hat{f}(x)=  \int_{\R^m} e^{-2\pi i<x,y>} f(y)dy. $$
	Using this  i, one has $ \norm{\hat{f}}_{L^2(\R^m,\R_m)} = \norm{f}_{L^2(\R^m,\R_m)}$. The inversion formula is then written as \cite{delanghe2012clifford}$$ \displaystyle f(x) =  \int_{\R^m} e^{2\pi i<x,y>} \hat{f}(y)dy. $$
	\begin{lem}\label{Fourier_inv}
		The Fourier transform of Clifford PSWFs is given by 
		$$ \displaystyle \widehat{\ps}(x) = \frac{(-1)^k}{c^d \mu_{N,m}^{k,c} } \ps\big(\frac{x}{c}\big) \chi_{(\overline{B(1)})} \big(\frac{x}{c}\big).$$
	\end{lem}
	
	\begin{proof}
		By the inverse Fourier transform, one has for a continuous $ f\in L^2(\R^m,\R_m)$, 
		\begin{eqnarray}\label{inv_bis}
			f(x) &=&  \int_{\R^m} e^{2 \pi i<x,z>} \ff f(z) d\g{z} =  \int_{\R^m} e^{2\pi i<x,z>} \int_{\R^m} e^{-2\pi i<z,y>} f(y)dy d\g{z} \nonumber \\
			&=&  \int_{\R^m} \int_{\R^m} e^{-2\pi i<\g{z},(y-x)>} f(y) dy d\g{z}.
		\end{eqnarray}
		Since
		$$ \ps(x) = \frac{1}{\mu^{(m)}_N(c)} \int_{\overline{B(1)}}e^{-2 \pi ic<x,y>} \ps(y)dy, $$
		one gets 
		\begin{eqnarray}
			\ff \ps(x)&=&\frac{1}{\mu^{(m)}_N(c)}\int_{\R^m} \int_{\overline{B(1)}} e^{-2\pi ic<\g{z},y>}\ps(y)dy e^{-i<\g{z},x>} dx \nonumber \\
			&=& \frac{1}{\mu^{(m)}_N(c)}\int_{\R^m} \Bigg( \int_{\R^m} e^{-2 \pi ic<\g{z},y-x>}\ps(\g{-y}/c)\chi_{(\overline{B(1)})}(y/c)dy\Bigg)  dx .\nonumber \\
		\end{eqnarray}
		To conclude the proof, it suffices to use the last equation together with \eqref{inv_bis}.
	\end{proof}
	
	\begin{lem}
		Let $c>0$ and let $m,k$ be two positive integers such that $c^2>\left( k+\frac{m}{2}-1\right)^2 - \frac{1}{4}$. Then, there exists an integer $N_0$ such that for all $n\geq N_0$, one has
		\begin{equation}\label{max-psi}
			\norm{\ps}_{\infty} \lesssim \left[ \left( 2N+k+\frac{m-1}{2} \right)\left( 2N+k+\frac{m+1}{2} \right) \right]^{1/2}.    \end{equation}
	\end{lem}
	The previous lemma is a simple reproduction of Theorem 2 in \cite{BK2017} together with \ref{relation_circular_case}.
	
	\begin{thm}
		Let \(c>0\). Let \(f\) be a \(c\)-band-limited Clifford-valued function in \(L^2(B(1),\mathbb{C}_m)\).
		For integers \(N,M\ge 0\) define the partial CPSWF expansion
		\[
		S_N^{(M)} f
		:= \sum_{k=0}^{M}\sum_{n=0}^{N}\sum_{\ell=1}^{d_{k,m}}
		\langle f,\psi^{(k,c)}_{n,m,\ell}\rangle\,\psi^{(k,c)}_{n,m,\ell},
		\]
		where \(\{\psi^{(k,c)}_{n,m,\ell}\}\) denotes an orthonormal CPSWF basis and \(d_{k,m}=\dim \mathcal{M}^+_l(k)\).
		There exists a constant \(C_{m,c}>0\), depending only on \(m\) and \(c\), such that for all \(N,M\) with \(N,M>\frac{ec}{2}\) and every \(k\le M\),
		\begin{equation}\label{eq:thm4.3}
			\| f - S_N^{(M)} f \|_{L^2(B(1))}
			\le
			C_{m,c}\,(2N+k)\,
			\left(
			\frac{e c}{4N + 2k + m + 3}
			\right)^{\frac{2N+k+m}{2}} \|f\|_{L^2(B(1))}.
		\end{equation}
	\end{thm}
	
	\begin{proof}
		By Parseval's equality, 
		\begin{equation}\label{1}
			\displaystyle \norm{f-S^{(M)}_N.f}^2_2 = \sum_{m>M}\sum_{k>N}\sum_{\ell=1}^{N(m,k)} |<f,\ps>|^2.
		\end{equation}
		Then, it remains to estimate $|<f,\ps>|$.\\
		Recall that from the Fourier inversion formula, one has $$ \displaystyle f(\g{x}) = c^m \int_{\overline{B(1)}} e^{2\pi ic<\g{x},\g{y}>} \hat{f}(\g{y}) d\g{y}. $$
		Consequently, for any positive integer $k$, we have
		\begin{eqnarray}\label{2}
			<f,\ps> &=& \int_{\overline{B(1)}}f(\g{x}) \ps(\g{x}) d\g{x}\nonumber \\ 
			&=& c^m \int_{\overline{B(1)}} \ps(\g{x}) \int_{\overline{B(1)}} e^{2\pi ic<\g{x},\g{y}>} \hat{f}(\g{y}) d\g{y} d\g{x} \nonumber \\
			&=& c^m\int_{\overline{B(1)}}\hat{f}(\g{y}) \int_{\overline{B(1)}} e^{2\pi ic<\g{x},\g{y}>}\ps(\g{x})d\g{x} d\g{y} \nonumber \\
			&=&c^m\left|\mu^{(m)}_N(c)\right| \int_{\overline{B(1)}}\hat{f}(\g{y})\ps(\g{y}) d\g{y} \nonumber \\
			&\leq& c^m \frac{\pi^{m/4}}{\sqrt{\Gamma(m/2+1)}} \left|\mu^{(m)}_N(c)\right| \norm{\ps(\g{y})}_\infty \norm{f}_{L^2(\R^m,\R_m)}.\nonumber \\
		\end{eqnarray}
	\end{proof}
	
	\section{Numerical Simulations}
	
	In this section, we illustrate the approximation properties of the Clifford
	prolate spheroidal wave functions (CPSWFs). All computations are performed
	on the unit disk $B(1)$ in dimension $m=2$. We compare the CPSWF
	approximation with the classical Fourier–Bessel expansion, which is a
	standard reference method for problems involving radial band-limited functions.
	\subsection{Example 1} 
\begin{figure}[H] 
		\centering
		\begin{subfigure}{0.4\textwidth}
			\centering
			\includegraphics[width=\linewidth]{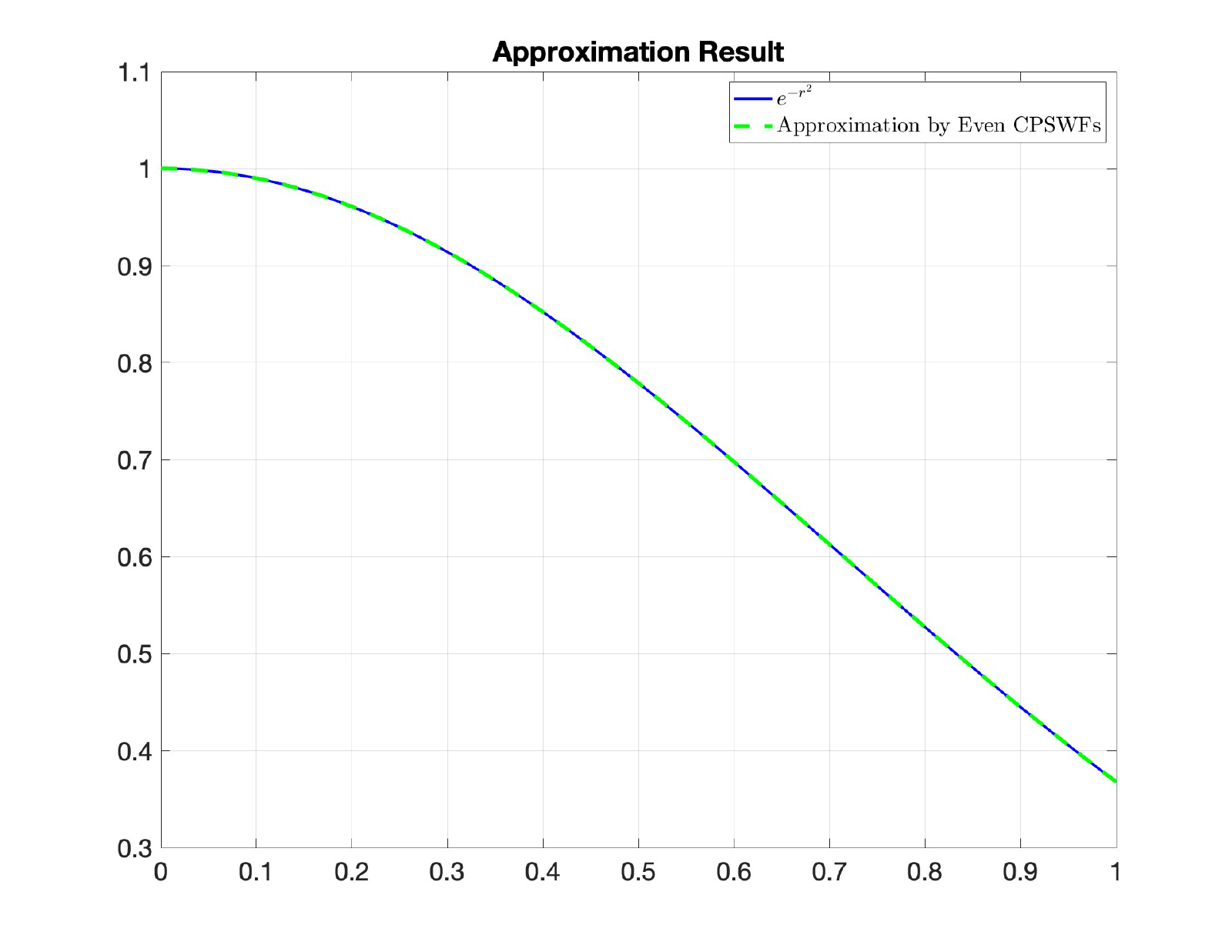} 
		\end{subfigure}
		\hfill    
		\begin{subfigure}{0.4\textwidth}
			\centering
			\includegraphics[width=\linewidth]{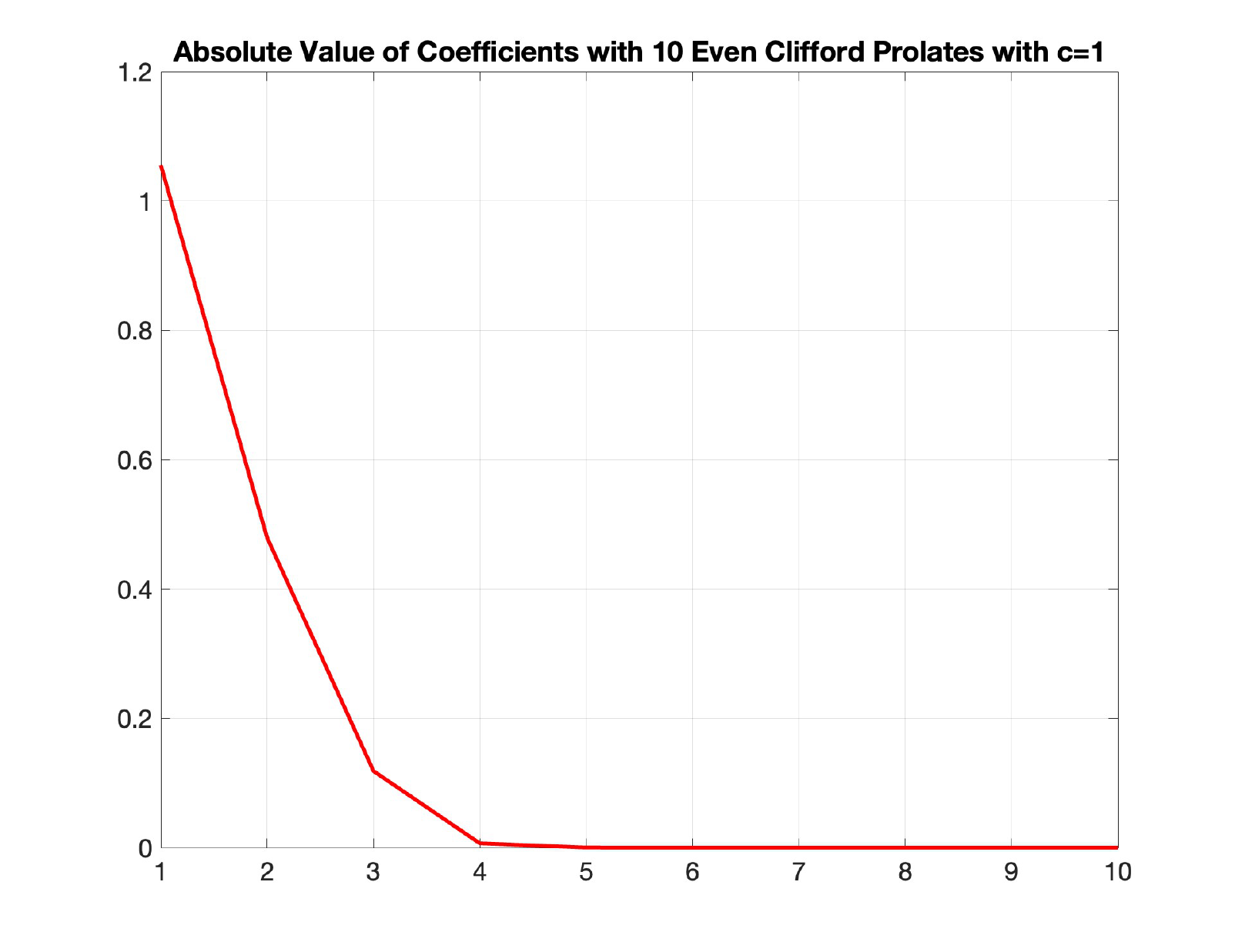} 
		\end{subfigure}
		\hfill 
		\begin{subfigure}{0.4\textwidth}
			\centering
			\includegraphics[width=\linewidth]{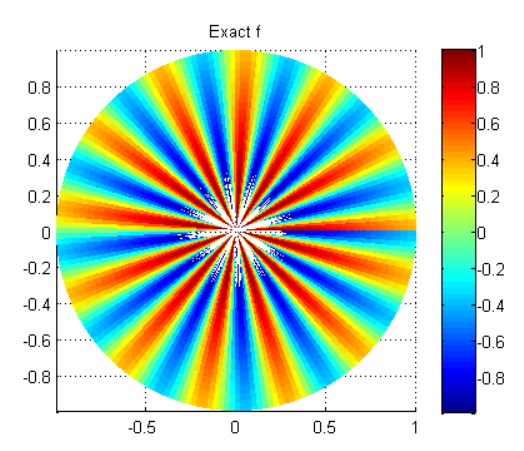} 
		\end{subfigure}
		\hfill
		\begin{subfigure}{0.40\textwidth}
			\centering
			\includegraphics[width=\linewidth]{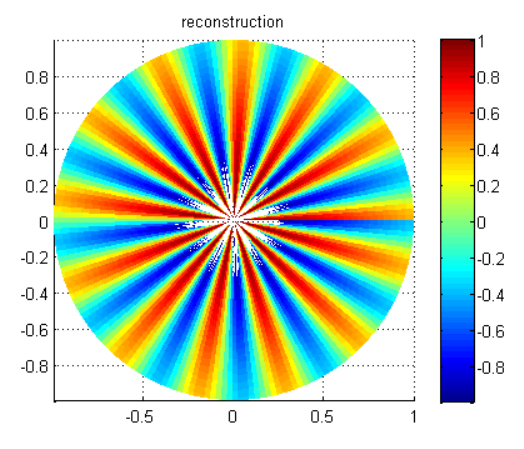} 
		\end{subfigure}
		\caption{The Graph of the exact function versus its representation using the $5$ first CPSWFs with $c=1$ and $n=0$ in Cartesian coordinates. } \label{example1}
	\end{figure}

In this example, we consider the function
	\[
	f(r,\theta) = e^{-r^{2}} \cos(4\pi\theta).
	\]
	We first reconstruct $f$ from its expansion coefficients in the basis of Clifford prolate spheroidal wave functions (CPSWFs). We then compare this reconstruction with the one obtained from the expansion of $f$ in the Fourier--Bessel basis as in \cite{Fourier_Bessel}. 
	Figure~\ref{example1} displays the exact function together with its
	CPSWF reconstruction using the first five basis functions. The CPSWF approximation captures the oscillatory angular structure with high fidelity.

	Table \ref{Comparaison table} reports the $L^2$ approximation error as a function of the number of basis functions. The CPSWF approximation converges much faster than the Fourier–Bessel expansion: with only $5$ CPSWFs, the accuracy is already of order $10^{-6}$, whereas Fourier–Bessel series require many more terms to reach comparable precision.
	\begin{table}[h!] 
		\centering
		\begin{tabular}{c c c}
			\hline
			\textbf{Number of basis functions $n$} & \textbf{$L_2$ error (CPSWFs)} & \textbf{$L_2$ error (Fourier-Bessel)} \\
			\hline
			5  & $1.5\times10^{-6}$ & $0.13$ \\
			7 & $8.24\times10^{-12}$ & $0.06$ \\
			10 & $6.6\times10^{-12}$ & $0.05$ \\
			\hline
		\end{tabular}
		\caption{Comparison of the $L_2$ error as a function of the number of basis functions $n$ for two different bases.}
		\label{Comparaison table}
	\end{table}
	
	\begin{rem}
		Using $150$ Fourier–Bessel modes, the approximation error remains of order $10^{-3}$, highlighting the efficiency of the CPSWF basis for band-limited signals.
	\end{rem}
	
	\subsection{Example 2}
	We consider the function
	$$ F(r,\theta) = \frac{J_0(cr)}{1+r^2} + \frac{J_1(cr)}{1+r^2} \theta e_1=g(r)+h(r)\theta e_1,$$
	which combines radial Bessel behaviour with an angular monogenic component.
	The goal is to evaluate the quality of the CPSWF reconstruction for this
	non-scalar Clifford-valued function. \\
	For a 2D-radial function, it is well known that the Fourier transform is identified with the Hankel transform of order zero :
	$$ \ff[f](\xi)=\hat{f}(\xi) = 2\pi \int_0^\infty f(r)J_0(|\xi|r)r dr = 2\pi \mathcal{H}_0(f)(|\xi|) .$$
	We recall the following two key identities :
	\begin{itemize}
		\item  \begin{equation}\label{Four1}
			\displaystyle \ff[J_0(ar)](\xi) = 2\pi \frac{\delta(|\xi|-a)}{|\xi|}.
		\end{equation}
		\item \begin{equation}\label{Four2}
			\displaystyle \ff\left[\frac{1}{1+r^2}\right](\xi)=2\pi K_0(|\xi|).
		\end{equation}
	\end{itemize}
	Where $K_0$ is the modified Bessel function defined by $ \displaystyle K_0(z)=\int_0^\infty \frac{\cos(zt)}{\sqrt{1+t^2}} dt$ having the crucial property 
	\begin{equation}\label{Mod_Bessel_decay}
		K_0(z) \sim _{z\to \infty} \sqrt{\frac{\pi}{2z}}e^{-z}.
	\end{equation}
	
\begin{figure}[H]
		\centering
		\includegraphics[width=1.\linewidth]{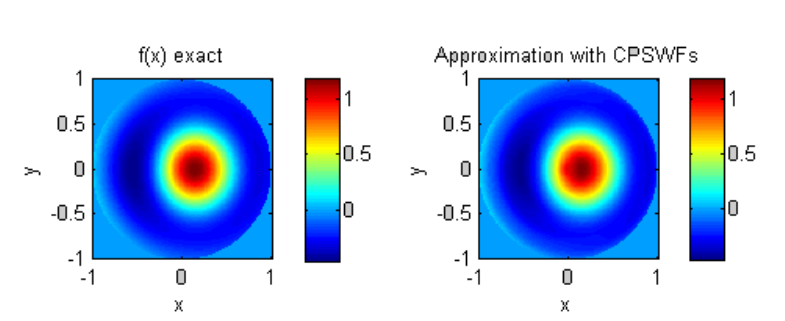}
		\caption{The graph of the exact $F$ versus its reconstruction using the 5 first CPSWFs.} \label{reconstruction4}
	\end{figure}

	Then, using \eqref{Four1} and \eqref{Four2}, one gets
	$$ \mathcal{F}[g](\xi) = 2\pi \left( \frac{\delta(|\xi|-a)}{|\xi|} \right) \ast \left( 2\pi K_0(|\xi|)\right) \approx \int_{|\omega|=c} K_0\left( |\xi-\omega|\right) d\sigma(\omega).$$
	Hence, using \eqref{Mod_Bessel_decay} and integrating over the circle
	$$ \mathcal{F}[g](\xi) \lesssim e^{-(|\xi|-c)}. $$
	Then, outside the disc $|\xi|\leq c$, $\hat{g}$ decay exponentially as $|\xi|$ increases. Note that these calculations are also valid for $h(r)$. We conclude that $F$ is almost band-limited. Figure~\ref{reconstruction4} illustrates the reconstruction of $F$ using the
	first five CPSWFs. The approximation captures both radial and angular
	components accurately.
	\FloatBarrier
	\section*{Acknowledgment}
	\noindent HBG was supported by a Lift-off fellowship from the Australian Mathematical Society and wishes to thank Dr Jeffrey Hogan for his helpful discussions and suggestions. \\
	\noindent AS was supported by the DGRST  research grant  LR21ES10.
	
	\bibliographystyle{elsarticle-num}
	\bibliography{Expansion_into_Clifford_Prolate_Spheroidal_Wave_Functions}
	
\end{document}